\documentclass{article}
\usepackage{amsmath,amsrefs,amsfonts,amssymb,amsthm}
\usepackage{accents}
\usepackage{authblk}
\usepackage[shortlabels]{enumitem}
\usepackage{geometry}
\usepackage{hyperref}
\usepackage{mathtools}
\usepackage{multicol}
\usepackage{tikz}
\usetikzlibrary{calc,decorations.pathreplacing}
\usepackage{pgffor}

\newtheorem{thm}{Theorem}[section]
\newtheorem*{thm*}{Theorem}
\newtheorem{lemma}[thm]{Lemma}
\newtheorem{conj}[thm]{Conjecture}
\newtheorem{cor}[thm]{Corollary}
\newtheorem{prop}[thm]{Proposition}

\theoremstyle{definition}
\newtheorem{defi}[thm]{Definition}

\theoremstyle{remark}
\newtheorem{convention}[thm]{Convention}
\newtheorem{example}[thm]{Example}
\newtheorem*{example*}{Example}
\newtheorem{remark}[thm]{Remark}

\newtheorem*{notation*}{Notation}

\newcommand{\area}{\mathsf{area}}
\newcommand{\dinv}{\mathsf{dinv}}
\newcommand{\revmaj}{\mathsf{revmaj}}
\newcommand{\matstat}{\mathsf{inv_3}}
\newcommand{\monot}{\mathsf{monot}}
\newcommand{\dw}{\mathsf{dw}}

\newcommand{\rev}{\mathsf{rev}}
\newcommand{\sched}{\mathsf{sched}}

\newcommand{\Des}{\mathsf{Des}}
\newcommand{\D}{\mathsf{D}} % Dyck paths
\newcommand{\LD}{\mathsf{LD}} % Labelled Dyck paths
\newcommand{\stLD}{\mathsf{stLD}} % standardly Labelled Dyck paths
\newcommand*{\dt}[1]{%
    \accentset{\mbox{\LARGE\bfseries .}}{#1}
    }
\newcommand{\Ht}{\widetilde{H}}

\newcommand{\N}{\mathbb{N}}

\title{Combinatorics of the Delta conjecture at $q=-1$}
\author[1]{Sylvie Corteel}
\author[1]{Matthieu Josuat-Vergès}
\author[2]{Anna Vanden Wyngaerd}
\affil[1]{IRIF, CNRS et Université Paris-Cité}
\affil[2]{Université Libre de Bruxelles}

\begin{document}
    \maketitle
    \abstract{
        In the context of the shuffle theorem, many classical integer sequences appear with a natural refinement by two statistics $q$ and $t$: for example the Catalan and Schröder numbers.  In particular, the bigraded Hilbert series of diagonal harmonics is a $q,t$-analog of $(n+1)^{n-1}$ (and can be written in terms of symmetric functions via the nabla operator).  The motivation for this work is the observation that at $q=-1$, this $q,t$-analog becomes a $t$-analog of Euler numbers, a famous integer sequence that counts alternating permutations.  We prove this observation via a more general statement, that involves the Delta operator on symmetric functions (on one side), and new combinatorial statistics on permutations involving peaks and valleys (on the other side).  An important tool are the schedule numbers of a parking function first introduced by Hicks; and expanded upon by Haglund and Sergel.  Other empirical observation suggest that nonnegativity at $q=-1$ holds in far greater generality.
    }        

\tableofcontents

\section{Introduction}
In the early 2000's, Haglund, Haiman, Remmel, Loehr and Ulyanov stated the \emph{shuffle conjecture}~\cite{HaglundHaimanLoehrRemmelUlyanov2005}: a combinatorial formula for the symmetric function $\nabla e_n$ in terms on \emph{labelled Dyck paths}. 
The interest in the symmetric function $\nabla e_n$ (where $\nabla$ is the MacDonald eigenoperator introduced in \cite{BergeronGarsia1998}) stems from it being the bi-graded Frobenius characteristic of the diagonal harmonic representation of the symmetric group~\cite{Haiman2002}. More than a decade after its statement, Carlsson and Mellit proved the full shuffle conjecture, which thus became a theorem~\cite{CarlssonMellit2018}. By then, many special cases were known: for example $\langle \nabla e_n, e_n\rangle$ gives the famous \emph{$q,t$-Catalan numbers} \cite{Haglund2003} and $\langle \nabla e_n, h_de_{n-d}\rangle$ the \emph{$q,t$-Schröder numbers}~\cite{Haglund2004}. A consequence of the full shuffle theorem is that the bi-graded Hilbert series $\langle \nabla e_n, h_1^n\rangle$ gives a $q,t$-analogue of $(n+1)^{n-1}$.  It can be described combinatorially as the generating function of length $n$ parking functions with respect to area and number of diagonal inversions.

The famous {\it Euler numbers} $(E_n)_{n\geq0}$ can be defined by their generating series:
\[
    \sum_{n\geq 0} E_n \frac{z^n}{n!}
        =
    \tan(z) + \sec(z).
\]
They answer various enumeration problems, the most famous one being that $E_n$ is the number of alternating permutations in $\mathfrak{S}_n$, that is, those $\sigma$ such that $\sigma_1 > \sigma_2 < \sigma_3 > \cdots$.  They also appear in Arnold's theory of singularity, and in number theory via their relation with Bernoulli numbers.  Seeing them in the context of Macdonald $q,t$-combinatorics is new, and is the motivation for this project.  We will show that specializing $q=-1$ in the $q,t$-analog of $(n+1)^{n-1}$ (the bi-graded Hilbert series of diagonal harmonics) gives:
\begin{align}
    \left. \langle \nabla e_n,h_1^n \rangle\right|_{q=-1} =  t^{\lfloor n^2/4\rfloor}E_n(t) \label{eq:t-euler}
\end{align}
where $E_n(t)$ is a $t$-analogue of $E_n$ appearing in \cite{HanRandrianarivonyZeng1999}.  This specialization at $q=-1$ is a $t$-refinement of the identity 
\begin{equation}
  \sum_{P \in \mathsf{PF}_n} (-1)^{\area(P)}
  =
  E_n,\label{eq:minus1-to-area}
\end{equation}
where $\mathsf{PF}_n$ are the parking functions of size $n$. The history of this identity can be found in \cites{MallowsRiordan1968,Kreweras1980,Pansiot1982}. For the definition of parking functions and their correspondence to standardly labelled Dyck paths, see \cite{Haglund2008}*{Chapter~5}.

We will establish Equation~\eqref{eq:t-euler} as a corollary of a more general statement involving a generalization of the shuffle theorem: the \emph{valley version} of the \emph{Delta conjecture}~\cite{HaglundRemmelWilson2018}. This is a combinatorial formula for the symmetric function $\Delta'_{e_{n-k-1}}e_n$. We will mainly use the following consequence of the Delta conjecture: 
\begin{equation}
    \langle \Delta'_{e_{n-k-1}}e_n,h_1^n \rangle= \sum_{P\in \stLD^{\bullet k}(n)}q^{\dinv(P)}t^{\area(P)} x^P, 
    \label{eq:hilbert-delta}
\end{equation}
where $\stLD(n)^{\bullet k}$ denotes the set of \emph{standardly labelled Dyck paths} with $k$ \emph{decorated valleys} and $\dinv$ and $\area$ are combinatorial statistics on this set. See Section~\ref{sec:paths} for the precise combinatorial definitions.  At $k=0$, we have $\Delta'_{e_{n-1}} e_n = \nabla e_n$, and the Delta conjecture reduces to the shuffle theorem. 

Specializations of the shuffle theorem and Delta conjecture at $q=0$ or $q=1$ have been extensively studied (see \cite{GarsiaHaglundRemmelYoo2019} and \cite{Romero2016}, respectively). To our knowledge, apart from \eqref{eq:minus1-to-area}, nothing much was known about the specialization at $q=-1$.

We were inspired by the following remarkable symmetric function identity, which first appeared in \cite{DAdderioIraciVandenWyngaerd2020}*{Theorem~4.11, case $m=0$}
\begin{equation}
    \sum_{k = 0}^{n-1} (-q)^{k}\Delta'_{e_{n-k-1}} e_n = \left.\nabla e_n \right|_{q=0}.
\end{equation}
Taking the scalar product with $h_1^n$ and evaluating at $q=-1$, we obtain
\begin{equation}\label{eq:sumdeltas}
    \sum_{k=0}^{n-1}\left. \langle \Delta'_{e_{n-k-1}}e_n,h_1^n\rangle\right|_{q=-1} = \langle \left. \nabla e_n\right|_{q=0}, h_1^n \rangle = [n]_t!,
\end{equation}
where the second equality is an easy consequence of the shuffle theorem. 

Our main result is a combinatorial interpretation of the terms of this sum, conditional on~\eqref{eq:hilbert-delta}. 
\begin{thm}\label{thm:main}
    For all $n\in \N$, we have
    \[\sum_{k = 0}^{n-1}\left(\sum_{P\in \stLD(n)^{\bullet k}}t^{\area(P)}(-1)^{\dinv(P)} \right)z^k= \sum_{\sigma \in \mathfrak S_n} t^{\matstat(\sigma)}z^{\monot(\sigma)},\] where $\matstat$ is a new statistic on permutations generalising Chebikin's notion of \emph{alternating descents} \cite{Chebikin2008} and $\monot(\sigma)$ is the number of double ascents or descents of $\sigma$ (see Section~\ref{sec:permutations} for the precise definitions).
\end{thm}
Thus if the Delta conjecture is proven to be true, we will have the following symmetric function interpretation. 
\begin{cor}\label{cor:main}
    If Equation~\eqref{eq:hilbert-delta} is true, then for all $n\in \N$ we have:
    \begin{align*}
        \sum_{k=0}^{n-1}\left. \langle \Delta'_{e_{n-k-1}}e_n,h_1^n\rangle\right|_{q=-1}z^k = \sum_{\sigma\in\mathfrak S_n} t^{\matstat(\sigma)}z^{\monot(\sigma)}.
    \end{align*}
\end{cor}

Notice that at $z=1$ our theorem agrees with Equation~\eqref{eq:sumdeltas}. 

The specialisation at $z=0$ of our theorem, and the fact that at $k=0$ the Delta conjecture reduces to the shuffle theorem will imply our formula \eqref{eq:t-euler}.

Our proof relies on the schedule formula decomposition of the combinatorial side of the valley Delta conjecture provided in \cite{HaglundSergel2021}. We use this schedule framework to identify the valley decorated Dyck paths that do not cancel out when specializing to $q=-1$. We then provide a bijection between these paths and permutations. This map will be defined via specific generating trees of the objects and will send $\area$ to $\matstat$ and the number of decorations to $\monot$. In this way the paths with no decorations ($k=0$) get sent to permutations with no double ascents or descents, that is, alternating permutations.

\section{The valley Delta conjecture}\label{sec:valley-delta}
In this section, we give the definitions needed to state the valley Delta conjecture.

\subsection{Valley-decorated labelled Dyck paths}\label{sec:paths}
\begin{figure}[ht!]
    \centering
    \begin{tikzpicture}[scale = .6]
        \draw[draw=none, use as bounding box] (-1, -1) rectangle (9,9);
        \draw[step=1.0, gray!60, thin] (0,0) grid (8,8);

        \draw[gray!60, thin] (0,0) -- (8,8);

        \draw[blue!60, line width=1.5pt] (0,0) -- (0,1) -- (0,2) -- (0,3) -- (1,3) -- (2,3) -- (2,4) -- (3,4) -- (4,4) -- (4,5) -- (4,6) -- (5,6) -- (6,6) -- (6,7) -- (7,7) -- (7,8) -- (8,8);

        \node at (0.5,0.5) {$3$};
        \draw (0.5,0.5) circle (.4cm); 
        \node at (0.5,1.5) {$4$};
        \draw (0.5,1.5) circle (.4cm); 
        \node at (0.5,2.5) {$5$};
        \draw (0.5,2.5) circle (.4cm); 
        \node at (2.5,3.5) {$1$};
        \draw (2.5,3.5) circle (.4cm); 
        \node at (4.5,4.5) {$6$};
        \draw (4.5,4.5) circle (.4cm); 
        \node at (4.5,5.5) {$7$};
        \draw (4.5,5.5) circle (.4cm); 
        \node at (6.5,6.5) {$2$};
        \draw (6.5,6.5) circle (.4cm); 
        \node at (7.5,7.5) {$8$};
        \draw (7.5,7.5) circle (.4cm); 
        \node at (3-1-0.5,3+0.5) {$\bullet$};
        \node at (4-0-0.5,4+0.5) {$\bullet$};
        \node at (6-0-0.5,6+0.5) {$\bullet$};
        \node at (7-0-0.5,7+0.5) {$\bullet$};

    \end{tikzpicture}
    \caption{An element of $\stLD^{\bullet 4}(8)$.}
    \label{fig:dyckpath}
\end{figure}
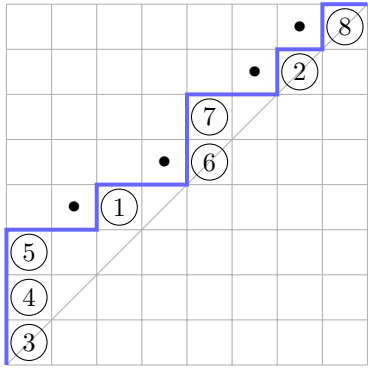

\begin{defi}
    A \emph{Dyck path} of size $n$ is a lattice path going from $(0,0)$ to $(n,n)$ consisting of east or north unit steps, always staying weakly above the line $x=y$, called the \emph{main diagonal}. The set of dyck paths is denoted by $\D(n)$.
\end{defi}

\begin{defi}
    A \emph{labelled Dyck path} is a pair $(\pi, w)$, where $\pi \in \D(n)$ and $w$ its \emph{labelling}: a word of positive integers whose $i$-th letter labels the $i$-th vertical step of $\pi$, placed in the square to the right of this step, such that the labels appearing in the same column are increasing from bottom to top. A labelling is said to be \emph{standard} if its the labels are exactly $1,2,\dots, n$. The set of (standardly) labelled Dyck paths of size $n$ is denoted by $\LD(n)$ (respectively, $\stLD(n)$).
\end{defi}
Standardly labelled Dyck paths are in bijection with parking functions. 

\begin{defi}
    The \emph{area word} of a Dyck path $\pi \in \D(n)$ is the word $a$ of $n$ non-negative integers whose $i$-th letter is the number of whole squares between the $i$-th vertical step of $\pi$ and the main diagonal $x=y$. The \emph{area} of a Dyck path is the sum of the letters of its area word and is denoted by $\area(\pi)$. 
\end{defi}

\begin{defi}
    Given $P\coloneqq (\pi,w)\in \LD(n)$ with area word $a$, the $i$-th vertical step of $P$ is called a \emph{contractible valley} if 
    \begin{itemize}
        \item either $a_{i-1}> a_{i}$,
        \item or $a_{i-1} = a_{i}$ and $w_{i-1} < w_{i}$. 
    \end{itemize}
%    In other words, it is either a vertical step preceded by two horizontal steps, or  a horizontal step preceded by a vertical step and if we were to delete (or ``contract'') the horizontal step preceding it, the labelling would remain valid. 
In other words, the $i$-th vertical step is a contractible valley if it is preceded by a horizontal step and the following holds:  after replacing the two steps
\tikz{\draw[blue!60, line width=1.5pt] (0,0)-|(.3,.3);} with \tikz{\draw[blue!60, line width=1.5pt] (0,0)|-(.3,.3);} (and accordingly shifting the $i$-th label one cell to the left), we still get a valid labelled path where labels are increasing in each column.
\end{defi}

\begin{defi}
    A \emph{(valley) decorated labelled Dyck path} is a triple $(\pi,w,dv)$ where $(\pi,w)\in \LD(n)$ and $dv$ is some subset of the contractible valleys of $(\pi,w)$. The elements of $dv$ are called \emph{decorations} and we visualise them by drawing a $\bullet$ to the left of these contractible valleys. The set $\LD(n)^{\bullet k}$ denotes the decorated labelled Dyck paths with exactly $k$ decorations. 
\end{defi}
\begin{defi}
    Given $P\coloneqq (\pi,w, dv) \in \LD(n)^{\bullet k}$ with area word $a$, a pair $(i,j)$ of indices of vertical steps with $1 < j\leq n$ is said to be a 
    \begin{itemize}
        \item \emph{primary diagonal inversion} if $a_i = a_j, w_i<w_j$ and  $i\not \in dv$,
        \item \emph{secondary diagonal inversion} if $a_i = a_j+1, w_i>w_j$ and  $i\not \in dv$.
    \end{itemize}
    The \emph{dinv} of $P$ is defined to be the total number of primary and secondary dinv pairs minus the number of decorated valleys and is denoted by $\dinv(P)$. 
\end{defi}

\begin{remark}
    We note that the $\dinv$ of a decorated labelled path is always a non-negative integer. Indeed, upon some reflection, one notices that each contractible valley forces the existence of at least one dinv pair.
\end{remark}

\begin{defi}
    Given $P\coloneqq (\pi,w, dv) \in \LD(n)^{\bullet k}$, the \emph{area} of $P$ is simply defined as the area of the underlying Dyck path, disregarding the labels and decorations: $\area(P)\coloneqq \area(\pi)$. 
\end{defi}

\begin{example}
    See Figure~\ref*{fig:dyckpath} for an example of an element of $\stLD(8)^{\bullet 4}$. Its labelling is $34516728$, its area word $01210100$ and so its area is $5$. 
    Its primary dinv pairs are 
    \[
        (1,5), (1,8), (2,6),
    \] 
    and its secondary dinv pairs are 
    \[
        (2,7),(3,4),(6,7).
    \] 
    Thus, since there are $4$ decorated valleys, the $\dinv$ is equal to $2$. 
\end{example}

\subsection{Symmetric functions}\label{sec:symfun}
For all the undefined notations and the unproven identities, we refer to \cite{DAdderioIraciVandenWyngaerd2022}, where definitions, proofs and/or references can be found. 

We denote by $\Lambda$ the graded algebra of symmetric functions with coefficients in $\mathbb{Q}(q,t)$, and by $\langle\, , \rangle$ the \emph{Hall scalar product} on $\Lambda$, defined by declaring that the Schur functions form an orthonormal basis. The standard bases of the symmetric functions are the monomial $\{m_\lambda\}_{\lambda}$, complete $\{h_{\lambda}\}_{\lambda}$, elementary $\{e_{\lambda}\}_{\lambda}$, power $\{p_{\lambda}\}_{\lambda}$ and Schur $\{s_{\lambda}\}_{\lambda}$ bases.

For a partition $\mu \vdash n$, we denote by \[ \tilde H_\mu \coloneqq \tilde H_\mu[X] = \tilde H_\mu[X; q,t] = \sum_{\lambda \vdash n} \widetilde{K}_{\lambda \mu}(q,t) s_{\lambda} \] the \emph{(modified) Macdonald polynomials}, where \[ \widetilde{K}_{\lambda \mu} \coloneqq \widetilde{K}_{\lambda \mu}(q,t) = K_{\lambda \mu}(q,1/t) t^{n(\mu)} \] are the \emph{(modified) Kostka coefficients} (see \cite{Haglund2008} for more details). 

Macdonald polynomials form a basis of the ring of symmetric functions $\Lambda$. This is a modification of the basis introduced by Macdonald \cite{Macdonald1998}.

\begin{figure}
	\centering
	\begin{tikzpicture}[scale=.5]
		\draw[gray] (0,0) grid (6,1);
		\draw[gray] (0,1) grid (6,2);
		\draw[gray] (0,2) grid (4,3);
		\draw[gray] (0,3) grid (4,4);
		\draw[gray] (0,4) grid (4,5);
		\draw[gray] (0,5) grid (3,6);
		\node at (2.5,3.5) {$c$};
		\fill[blue, opacity=.4] (0,3) rectangle (2,4) node[midway, opacity=1, black]{co-arm};
		%\fill[blue, opacity=.4] (3,3) rectangle (7,4) node[midway, opacity=1, black]{arm};
		%\fill[blue, opacity=.4] (2,4) rectangle (3,7) node[midway, opacity=1, black, rotate=90]{leg};
		\fill[blue, opacity=.4] (2,3) rectangle (3,0) node[midway, opacity=1, black, rotate=90]{co-leg};
	\end{tikzpicture}
	\caption{co-arm and co-leg of a cell in a partition.}\label{fig:limbs}
\end{figure}
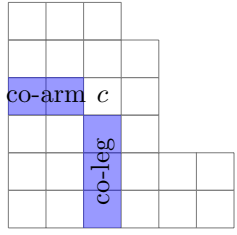

If we identify the partition $\mu$ with its Ferrer diagram, i.e.~with the collection of cells $\{(i,j)\mid 1\leq i\leq \mu_j, 1\leq j\leq \ell(\mu)\}$, then for each cell $c\in \mu$ we define the \emph{co-arm} and \emph{co-leg} (denoted respectively as $a'_\mu(c), l'_\mu(c)$) as the number of cells in $\mu$ that are strictly to the left and below $c$ in $\mu$, respectively (see Figure~\ref{fig:limbs}).  
Define the following constant:
\begin{align*}
	B_{\mu} & \coloneqq B_{\mu}(q,t) = \sum_{c \in \mu} q^{a_{\mu}'(c)} t^{l_{\mu}'(c)}.
\end{align*}
Let $f[g]$ denotes the \emph{plethystic evaluation} of a symmetric function $f$ in an expression $g$  (see \cite{Haglund2008}*{Chapter~1 page 19}).
\begin{defi}[{\cite{BergeronGarsia1998}*{3.11}}]
	\label{def:nabla}
	We define the linear operator $\nabla \colon \Lambda \rightarrow \Lambda$ on the eigenbasis of Macdonald polynomials as \[ \nabla \Ht_\mu = e_{\lvert \mu \rvert}[B_\mu] \Ht_\mu=q^{n(\mu')}t^{n(\mu)} \Ht_\mu.\]
\end{defi}

\begin{defi}
	\label{def:delta}
	For $f \in \Lambda$, we define the linear operators $\Delta_f, \Delta'_f \colon \Lambda \rightarrow \Lambda$ on the eigenbasis of Macdonald polynomials as \[ \Delta_f \Ht_\mu = f[B_\mu] \Ht_\mu, \qquad \qquad \Delta'_f \Ht_\mu = f[B_\mu-1] \Ht_\mu. \]
\end{defi}

\subsection{The statement}

Now we have all the necessary definitions to state the valley Delta conjecture, which first appeared in \cite{HaglundRemmelWilson2018}.
\begin{conj}\label{conj:Delta}
    For all $n, k\in \N$,
    \begin{align*}
        \Delta'_{e_{n-k-1}}e_n = \sum_{P\in \LD^{\bullet k}(n)}q^{\dinv(P)}t^{\area(P)} x^D
    \end{align*}
    where the sum is over the set of \emph{labelled Dyck paths} of size $n$ with $k$ decorations on \emph{contractible valleys}.
\end{conj}
Taking the (Hall) scalar product with $h_1^n$ (i.e.,~the Hilbert series), the Delta conjecture implies~\eqref{eq:hilbert-delta}.

\section{Schedule formula}
In this section we discuss the schedule formula for the combinatorics of the valley Delta conjecture proved by Haglund and Sergel in \cite{HaglundSergel2021}. Their formula is an extension of the first work on schedule numbers by Hicks in her thesis \cite{Hicks2013}. 

\begin{defi}
    The set of \emph{decorated permutations} $\mathfrak S_n^{\bullet k}$ of $[n] := \{1,\dots,n\}$ is the set of permutations of $[n]$ where $k$ of its $n$ letters are \emph{decorated}, represented as $\dt \sigma_i$. Set $\mathfrak S_n^\bullet \coloneqq \sqcup_{k\in [n]} \mathfrak S_n^{\bullet k}$.  For $\sigma \in \mathfrak{S}_n^{\bullet}$, we denote by $\mathsf{dec}(\sigma)$ its number of decorations. 
\end{defi}
\begin{defi}
    Given $P\coloneqq (\pi, w,dv)\in \stLD(n)^{\bullet k}$ with area word $a$, the \emph{diagonal word} of a decorated labelled Dyck path is the decorated permutation of $n$ obtained as follows. A label $w_i$ of $P$ is said to lie in the \emph{$j$-th diagonal} if $a_i = j$. List all the labels $w_i$ in the $0$-th diagonal, in decreasing order, adding a decoration on the label if $i\in dv$. Then do the same for the $1$-st diagonal, $2$-nd diagonal, and so forth. Denote this diagonal word by $\dw(P)$. 
\end{defi}

\begin{example}
    The diagonal word of the path in Figure~\ref{fig:dyckpath} is $\dt 8\dt 63\dt 2 74\dt 15$. 
\end{example}

\begin{defi}\label{def:maj}
    For $\sigma\in \mathfrak S_n$ a permutation, its \emph{major index} is defined to be the sum of the elements of the set $\{i\in [n-1]\mid \sigma(i)>\sigma(i+1)\}$. The \emph{reverse major index} of a permutation is simply the major index of the reverse permutation $\sigma^{\mathsf{rev}}\coloneqq \sigma_{n}\cdots \sigma_1$. For any marked permutation $\sigma\in \mathfrak S_n^\bullet$, denote by $\revmaj(\sigma)$ the reverse major index of its underlying permutation. 
\end{defi}

The following is an easy consequence of the definitions. 
\begin{prop}
    For all $P\in \stLD(n)^{\bullet k}$, we have $\area(P) = \revmaj(\dw(P))$. 
\end{prop}
\begin{proof}
    Due to the condition of strictly increasing labels in the columns of labelled Dyck paths, each diagonal has at least one label which is bigger than some label of the previous diagonal. Thus, the labels in the $j$-th diagonal of $P$ are exactly the numbers in the $(j+1)$-th decreasing run of $\dw(P)$, and they each contribute $j$ units to the area. 
\end{proof}

The following convention will greatly simplify definitions and proofs.
\begin{convention}\label{conv:preceding-zero}
    Given $\sigma\in \mathfrak{S}_n$, we will implicitly consider that $\sigma$ is preceded by a $0$-th entry: $\sigma_0=0$. If $\sigma\in \mathfrak{S}_n^\bullet$, then $\sigma_0$ is never decorated.  
\end{convention}

\begin{defi}[Hicks~\cite{Hicks2013}]\label{def:schedule}
    For $\tau\in\mathfrak S_n^\bullet$  define its \emph{schedule numbers} $\sched(\tau)= (s_i)_{1\leq i \leq n}$ as follows. Take $r_0, r_1,r_2,\dots$ to be the decreasing runs of $\tau_0\cdots \tau_n $ (we have $r_0 = 0$). 
    \begin{itemize}
        \item If $\tau_i$ is undecorated and an element of $r_j$, let 
        \begin{align*}
            s_i = &\#\{k\in r_j\mid k \text{ is undecorated and } k>\tau_i\} 
            \\ &+ \#\{k\in r_{j-1} \mid k \text{ is undecorated and } k<\tau_i\}.
        \end{align*}
        \item If $\tau_i$ is decorated and an element of $r_j$, let  
        \begin{align*}
            s_i =& \#\{k\in r_j\mid k \text{ is undecorated and } k<\tau_i\} 
            \\ & + \#\{k\in r_{j+1} \mid k \text{ is undecorated and } k>\tau_i\}.
        \end{align*}
    \end{itemize}
\end{defi}

To reformulate this definition of schedule numbers, we introduce the following.

\begin{defi} \label{def:cyclic-runs}
	Let $\sigma \in \mathfrak{S}_n$.  A sequence of consecutive elements $\sigma_i,\dots,\sigma_j$ (with $0 \leq i \leq j$) in $\sigma$ is a {\it cyclic (decreasing) run} if there exists an integer $k$ such that $\sigma_i+k \mod n+1,\dots, \sigma_j+k\mod n+1$ is decreasing (where the modulo means we take the representative in $\{0,\dots,n\}$).
	Moreover, a cyclic run $\sigma_i,\dots,\sigma_j$ is {\it left-maximal} if either $i=0$ or $\sigma_{i-1},\dots,\sigma_j$ is not a cyclic run, and {\it right-maximal} if either $j=n$ or $\sigma_i,\dots,\sigma_{j+1}$ is not a cyclic run.
\end{defi}

\begin{example}
    Let $\sigma = 0649751832$.  Some right-maximal cyclic runs are $064$, $6497$, $7518$.  Some left-maximal cyclic runs are $1832$, $6497$.  
\end{example}
Note that for each $j$, there is a unique left-maximal cyclic run $\sigma_i,\dots,\sigma_j$ (obtained by choosing $i$ to be minimal), and similarly for each $i$ there is a unique right-maximal cyclic run $\sigma_i,\dots,\sigma_j$ (obtained by choosing $j$ to be maximal).

\begin{remark}\label{rem:cyclic-runs}
	The definition of schedules is rephrased as follows: if $\tau_i$ is undecorated (respectively, decorated), then $s_i$ is the number of undecorated values (excluding $\tau_i$) in the maximal decreasing cyclic run ending (respectively, starting) at $\tau_i$.  
\end{remark}

\begin{example}\label{ex:schedules}
    If $\tau$ is the diagonal word of the path in Figure~\ref*{fig:dyckpath}, then we have
    \begin{center}
        \begin{tabular}{lllll}
        $\tau$ &$0$ &$\dt 8\,\dt 6\,3\,\dt 2\,$ &$7\,4\,\dt 1\,$ &$5$\\
        $s$ &&$1\,2\,1\,2\,$ &$1\,2\,1\,$ & $1$.
        \end{tabular}
    \end{center}
    In Figure~\ref{fig:cyclic-runs}, we visualise $\tau$ by placing dots at coordinates $(i,\tau_i)$: white dots for undecorated and black dots for decorated $\tau_i$. Notice that if we view this picture as a cylinder, identifying the top and the bottom, the operation $\sigma_i\mapsto \sigma_i + k \mod n$ can be seen as a rotation of this cylinder, hence the name ``cyclic runs''. The maximal cyclic run starting at $\dt 6$ is $\dt 6 3 \dt 2 7$ (Figure~\ref{fig:cyclic-runs}, left), so that the schedule number is $2$. The maximal cyclic run ending at $4$ is $3\dt 274$ (Figure~\ref{fig:cyclic-runs}, right), so the schedule number is also $2$. 

\end{example}
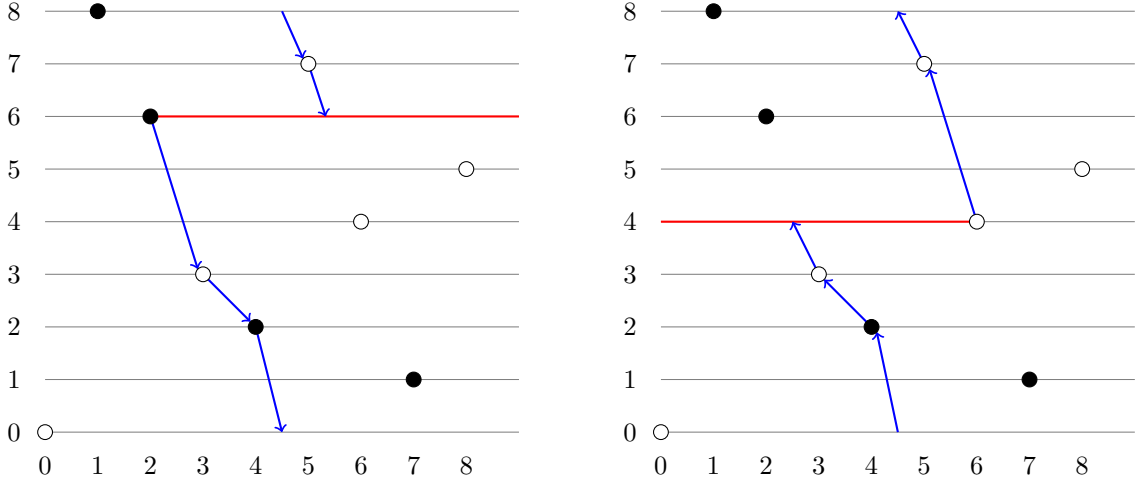
\begin{figure}
    \centering
    \begin{tikzpicture}[scale=.7]
        \foreach \i in {0,...,8}{ 
        \draw[opacity=.5] (0,\i) -- (9,\i);
        \node[left=.2cm] at (0,\i) {\i};
        \node[below=.2cm] at (\i,0) {\i};
        }        
        \draw[thick,red](2,6) -- (9,6);
        \draw[->, thick, blue] (2,6) -- (2.9,3.1);
        \draw[->, thick, blue] (2.9,3.1)--(3.9,2.1);
        \draw[->, thick, blue] (4,2)--(4.5,0);
        \draw[->, thick, blue] (4.5,8)--(4.9,7.1);
        \draw[->, thick, blue] (5,7)--(5.33,6);

        \draw[fill = white] (0,0) circle (4pt);
        \draw[fill] (1,8) circle (4pt);
        \draw[fill] (2,6) circle (4pt);
        \draw[fill = white] (3,3) circle (4pt);
        \draw[fill] (4,2) circle (4pt);
        \draw[fill = white] (5,7) circle (4pt);
        \draw[fill = white] (6,4) circle (4pt);
        \draw[fill] (7,1) circle (4pt);
        \draw[fill = white] (8,5) circle (4pt);
    \end{tikzpicture}
    \hfill
    \begin{tikzpicture}[scale=.7]
        \foreach \i in {0,...,8}{ 
        \draw[opacity=.5] (0,\i) -- (9,\i);
        \node[left=.2cm] at (0,\i) {\i};
        \node[below=.2cm] at (\i,0) {\i};
        }        
        \draw[thick,red](6,4) -- (0,4);
        \draw[->,thick,blue] (6,4)--(5.1, 6.9); 
        \draw[->,thick,blue] (5,7)--(4.5, 8);
        \draw[->,thick,blue] (4.5,0)--(4.1, 1.9); 
        \draw[->,thick,blue] (4,2)--(3.1, 2.9);
        \draw[->,thick,blue] (3,3)--(2.5,4); 
        \draw[fill = white] (0,0) circle (4pt);
        \draw[fill] (1,8) circle (4pt);
        \draw[fill] (2,6) circle (4pt);
        \draw[fill = white] (3,3) circle (4pt);
        \draw[fill] (4,2) circle (4pt);
        \draw[fill = white] (5,7) circle (4pt);
        \draw[fill = white] (6,4) circle (4pt);
        \draw[fill] (7,1) circle (4pt);
        \draw[fill = white] (8,5) circle (4pt);
    \end{tikzpicture}
    \caption[]{Illustration of cyclic runs.}\label{fig:cyclic-runs}
\end{figure}

The following result can be deduced from \cite{HaglundSergel2021}~{Theorem~3.13}. 
\begin{thm}\label{thm:schedule-formula}
    Given $\tau\in\mathfrak S_n^{\bullet k}$ and $(s_i)_{1\leq i \leq n}$ its schedule numbers, we have:
    \[
      \sum_{\substack{P\in \stLD(n)^{\bullet k}\\ \dw(P) = \tau}} q^{\dinv(P)}t^{\area(P)}= t^{\revmaj(\tau)} \prod_{i\in [n]} [s_i]_q.
    \]
\end{thm}

A consequence of this theorem is that a marked permutation is the diagonal word of some permutation if and only if all its schedule numbers are strictly positive.  At $q=-1$, more terms vanish as a consequence of the following.

\begin{lemma}\label{lem:schedules-are-interval}
	Let $\sigma\in \mathfrak S_n^\bullet$, and $(s_1,\dots, s_n) = \sched(\sigma)$.  If $s_i>0$ for all $i\in [n]$, then there exists $j\geq 1$ such that the set of schedule numbers $\{ s_i \; | \; 1\leq i \leq n \}$ is equal to $[j]$.
\end{lemma}

\begin{proof}
	Suppose there is an index $i$ such that $s_i>1$.  Our goal is to show that there exists $i'$ such that $s_{i'} = s_i-1$. We distinguish two cases, (which are very similar):
\begin{itemize}
	\item $\sigma_i$ is not decorated.  Consider the left-maximal run ending at $\sigma_i$, denoted $R$, so $R = (\sigma_h,\dots,\sigma_i)$ for some $h$.  Consider the maximal $i'$ such that $\sigma_{i'}$ is not decorated, and $h \leq i'<i$.  It exists because $s_i>0$.  Let $R' = \sigma_{h'},\dots,\sigma_{i'}$ be the left-maximal cyclic run ending at $\sigma_{i'}$.  We have $h'\leq h$ (because we know that $(\sigma_h,\dots,\sigma_{i'})$ is a cyclic run), so $R$ and $R'$ overlap on every non-decorated value of $R$ other than $\sigma_i$.  We deduce that $s_{i'} \geq s_i-1$.  In case of equality, we are done.  Otherwise, we iterate this construction to find $i''<i'$ such that $s_{i''} \geq s_{i'}-1$, etc.  Since the increment is at most 1, after some steps we find $i^{(k)}$ such that $s_{i^{(k)}} = s_{i}-1$.

	\item $\sigma_i$ is decorated.  Consider the right-maximal run beginning at $\sigma_i$, denoted $R$, so $R = (\sigma_i,\dots,\sigma_j)$ for some $j$. Consider the maximal $i'$ such that $\sigma_{i'}$ is not decorated, and $i < i' \leq j$.  It exists because $s_i>0$.  Let $R' = \sigma_{h'},\dots,\sigma_{i'}$ be the left-maximal cyclic run ending at $\sigma_{i'}$.  We have $h' \leq i$ (because we know that $(\sigma_{i},\dots,\sigma_{i'})$ is a cyclic run), so $R$ and $R'$ overlap on every non-decorated value of $R$ other than $\sigma_{i'}$.  We deduce that $s_{i'} \geq s_i-1$.  In case of equality, we are done.  Otherwise, we are back to the situation in the previous case. 
\end{itemize}
In either case, we eventually find $i'$ as announced.  The statement in the lemma follow by iteration.
\end{proof}

A consequence of this lemma is that the product $\prod_{i=1}^n [s_i]_q $ at $q=-1$ is $0$ unless all schedule numbers are $1$ (otherwise, one of the factors is $[2]_q=1+q$).  Therefore, the schedule formula of Theorem~\ref{thm:schedule-formula} becomes
\begin{equation}\label{eq:sched-q-1}
    \sum_{\substack{P\in \stLD(n)^{\bullet k}\\ \dw(P) = \tau}} (-1)^{\dinv(P)}t^{\area(P)}= t^{\revmaj(\tau)}.  
\end{equation}
\begin{notation*}
    We denote by  $\mathfrak S^{\bullet}_n(1^n)$ to be the subset of $\mathfrak{S}_n^{\bullet}$ of marked permutations with schedule $1^n$. 
\end{notation*}
Now, summing Equation~\eqref{eq:sched-q-1} over all possible permutations with $k$ decorations $\tau$, we get the following interpretation of the combinatorics of the Hilbert series of the Delta conjecture at $q=-1$:
\begin{equation}\label{eq:combinatorics-Delta-q-1}
    \sum_{P\in \stLD(n)^{\bullet k}} (-1)^{\dinv(P)}t^{\area(P)}= \sum_{\tau\in\mathfrak S_n^{\bullet k}(1^n)} t^{\revmaj(\tau)}.
\end{equation}

\begin{remark}
We announced in the introduction that in the case $k=0$ and $t=1$, the left-hand side of \eqref{eq:combinatorics-Delta-q-1} is the Euler number $E_n$.  Comparing with the right-hand side, this means that undecorated permutations $\sigma\in\mathfrak{S}_n$ with schedule $1^n$ are counted by the Euler number $E_n$.  On the other side, Ramassamy~\cite{Ramassamy2018}*{Corollary~5} gave a new combinatorial interpretation of the Euler number $E_n$ as the number of total cyclic orders on $\{0,\dots,n\}$ such that $(i,i+1,i+2)$ is clockwise oriented for each $i\in \{0,\dots,n-2\}$.  Undecorated permutations $\sigma\in\mathfrak{S}_n$ with schedule $1^n$ are simply related with Ramassamy's total cyclic orders via $\sigma \mapsto (0,\sigma^{-1}(1),\dots,\sigma^{-1}(n))$.
\end{remark}

We close this section by providing a bijection between $\mathfrak{S}_n$ and  $\mathfrak S^{\bullet}_n(1^n)$, whose inverse is given by simply removing the  decorations.  The existence of this bijection means that each permutation can be decorated in exactly one way, so that the result has schedule $1^n$.

\begin{lemma}
	Let $\sigma \in \mathfrak S^{\bullet}_n(1^n)$.  Let $\sigma_k,\dots,\sigma_\ell$ be a left-maximal cyclic run, and assume that $\sigma_\ell$ is undecorated.  Then $\sigma_k$ is undecorated as well, and $\sigma_{k+1},\dots,\sigma_{\ell-1}$ are decorated.
\end{lemma}

\begin{proof}
	Since $s_\ell=1$, there is exactly one undecorated entry in $\sigma_k,\dots,\sigma_{\ell-1}$.  Denote $\sigma_u$ this undecorated entry.  If $u\neq k$, $\sigma_k$ is decorated and has at least two decorated entries (namely $\sigma_u$ and $\sigma_\ell$) in the right-maximal cyclic run beginning at $\sigma_k$.  This would give $s_k\geq 2$, which contradicts $s_k=1$.  Thus $u=k$, and we get that $\sigma_k$ in the only non-decorated entry in $\sigma_k,\dots,\sigma_{\ell-1}$. 
\end{proof}

\begin{lemma}\label{lem:unique-decoration-1n}
	For each permutation $\sigma\in \mathfrak S_n$, there exists exactly one decorated permutation with underlying permutation $\sigma$ and schedule $1^n$.
\end{lemma}

\begin{proof}
	Assuming that there exists a decoration of $\sigma$ so that the schedule is $1^n$, the previous lemma readily gives necessary conditions on how to find it (in particular, uniqueness will follow from existence):  we proceed from right to left (starting from $\sigma_n$ and ending at $\sigma_1$), noting than $\sigma_n$ is necessarily undecorated (otherwise we have $s_n=0$).  We define a sequence of indices $i_1>i_2>\dots>i_m$ for some $m\geq 1$ as follows;
    \begin{itemize}
        \item $i_1 = n$,
        \item knowing $i_j$, we find $i_{j+1}$ by the condition that $\sigma_{i_{j+1}} \cdots \sigma_{i_j}$ is the left-maximal cyclic run ending at $\sigma_{i_j}$,
        \item the sequence stops at $i_m=0$.       
    \end{itemize}
    We claim that decorating the indices not in $\{i_1,\dots,i_m\}$ yields the unique decoration such that the associated schedule is $1^n$.  It remains only to check that the schedule of this decorated permutation is indeed $1^n$.

    By construction, we have $s_k = 1$ if $\sigma_k$ is undecorated.  It remains to show $s_k=1$ when $\sigma_k$ is decorated.  So, consider a right-maximal cyclic run $\sigma_k,\dots,\sigma_\ell$ where $\sigma_k$ is decorated.
    \begin{itemize}
	\item Suppose $s_k\geq 2$.  So, there are two (or more) undecorated entries in this run, say $\sigma_i$ and $\sigma_j$ with $k<i<j$.  The left-maximal cyclic run ending at $\sigma_j$ contains at least $\sigma_k,\dots,\sigma_j$, so $\sigma_i$ being undecorated contradicts the construction of $i_1>i_2>\dots>i_m$ as above.
        \item Suppose $s_k=0$. So, there are no undecorated entries in this run.  Let $\sigma_i$ be an undecorated entry with $i>k$ and $i$ minimal (this exists because $\sigma_n$ is undecorated).  We have $i>\ell$ (because $\sigma_{k+1},\dots,\sigma_\ell$ are decorated).  Consider the left-maximal cyclic run ending at $\sigma_i$, denoted  $\sigma_{i'}, \dots,\sigma_i$. It cannot begin at $\sigma_{i'}$ with $i' \leq k$, because $\sigma_k,\dots,\sigma_i$ is not a cyclic run (the right-maximal cyclic run beginning at $\sigma_k$ ends at $\sigma_\ell$, and $i>\ell$).  Thus $i'> k$, so $\sigma_{i'}, \dots,\sigma_i$ does not contain any undecorated entry apart from $\sigma_i$.  But this means $s_i=0$, which is a contradiction.
    \end{itemize}
    Other cases being excluded, we thus have $s_k=1$.  This completes the proof of existence and uniqueness of the decoration with schedule $1^n$.
\end{proof}

\section{Permutations}\label{sec:permutations}

We continue to use Convention~\ref{conv:preceding-zero}: for any $\sigma\in \mathfrak{S}_n$, we set $\sigma_0=0$. 

\begin{defi}
    Given a permutation $\sigma = \sigma_1\cdots \sigma_n$ and an index $i\in \{2,\dots,n\}$, we say that $\sigma_i$ is a 
    \begin{multicols}{2}
        \begin{itemize}
            \item \emph{double ascent} if $\sigma_{i-2} < \sigma_{i-1} < \sigma_i$;
            \item \emph{double descent} if $\sigma_{i-2} > \sigma_{i-1} > \sigma_i$;
            \item \emph{peak} if $\sigma_{i-2} < \sigma_{i-1} > \sigma_i$;
            \item \emph{valley} if $\sigma_{i-2} > \sigma_{i-1} < \sigma_i$. 
        \end{itemize}
    \end{multicols}
 
\end{defi}

\begin{defi}
    For $\sigma\in\mathfrak S_n$, we define 
    \[
        \monot(\sigma) = \#\big\{i\in \{2,\dots,n\} \; \big| \; \sigma_i \text{ is a double descent or a double ascent} \big\}.
    \]
\end{defi}

\begin{defi}\label{def:3inversion}
    Let $\sigma\in \mathfrak S_n$, a pair $(i,j)$ with $1\leq i<j\leq n$ is said to be a \emph{$3$-inversion} if one of the following holds:
    %TODO find better name for this?
    \begin{itemize}
        \item $\sigma_{j}$ is a double ascent and $\sigma_{j-1}<\sigma_i < \sigma_j$;
        \item $\sigma_{j}$ is a double descent and $\sigma_{j-1}> \sigma_i > \sigma_j$;
        \item $\sigma_{j}$ is a peak and $\sigma_i>\sigma_j$;
        \item $\sigma_{j}$ is a valley and $\sigma_i<\sigma_j$.
    \end{itemize}
    The number of $3$-inversions of $\sigma$ is denoted by $\matstat(\sigma)$. 
\end{defi}
\begin{example}
    For $n=3$, 123 has zero 3-inversions, 132 and 321 have one 3-inversions, 231 and 312 have two 3-inversion and  213 has three 3-inversions.  
\end{example} 
Though the definition of the statistic might not seem very natural, we will see in the proof of the proposition that it can be tracked via a rather simple insertion procedure on permutations (similar to the Lehmer codes).

We recall the following classical definitions.
\begin{defi}
    For $n\in \N$, define its \emph{$t$-analogue} by $[n]_t\coloneqq 1+t+\ldots +t^{n-1}$. The $t$-analogue of $n!$ is given by $[n]_t!=[n]_t[n-1]_t\ldots [1]_t$.
\end{defi}
\begin{defi}
    A statistics $I$ on permutation is called \emph{Mahonian} if
    $
    \sum_{\sigma \in S_n} t^{I(\sigma)}=[n]_t!.
    $
\end{defi}
Two classical Mahonian statistics are the major index (Definition~\ref{def:maj}) and the inversion number defined by
\[
    \mathsf{inv}(\sigma)
        =
    \#\big\{ (i,j) \; \big|\; 1\leq i<j\leq n \text{ and } \sigma_i>\sigma_j \big\}.
\]

From the main result in Section~\ref{sec:trees} (Theorem~\ref{thm:trees-isom}) and the fact that $\revmaj$ is Mahonian, we will be able to deduce the following. 
\begin{prop}\label{prop:matstat-mahonian}
The statistic $\matstat$ is Mahonian; that is
$$
\sum_{\sigma \in S_n} t^{\matstat(\sigma)}=[n]_t!.
$$
\end{prop}

 In \cite{Chebikin2008}, Chebikin defines another variant of the inversion statistics.
\begin{defi}
    Let $\sigma\in \mathfrak S_n$, define $\hat{c}_i(\sigma)$ to be the number of indices $j>i$ such that
    %TODO find better name for this?
    \begin{itemize}
        \item $i$ is odd and $\sigma_i>\sigma_j$; or
        \item $i$ is even and $\sigma_i<\sigma_j$.
    \end{itemize}
   Let $\hat{\imath}(\sigma)=\hat{c}_1(\sigma)+\ldots +\hat{c}_{n-1}(\sigma)$. 
\end{defi}

\begin{prop}[{\cite{Chebikin2008}*{Corollary 3.5}}]
The statistics $\hat{\imath}$ is Mahonian. Indeed we have
$$
\sum_{\sigma \in S_n} t^{\hat{\imath}(\sigma)}=[n]_t!.
$$
\end{prop}

\begin{defi}
    An \emph{alternating permutation} of $[n]$ is a permutation $\sigma \in \mathfrak{S}_n$ such that $\sigma_1>\sigma_2<\sigma_3>\cdots$. We denote the set of such permutations by $\mathfrak A_n$.
In other words, \[\mathfrak A_n \coloneqq \{\sigma \in \mathfrak S_n \mid \monot(\sigma) = 0\}.\]
\end{defi}

The alternating permutations are counted by the Euler numbers $E_n$, defined by the generating series 
\begin{equation*}
    \tan(x) + \sec(x) = \sum_{n\geq 0} E_n \frac{x^n}{n!}.
\end{equation*}

One can easily check that if $\sigma$ is alternating then $\hat{\imath}(\sigma)=\matstat(\sigma)$.
But this is not true in general. For example if $\sigma=123$, $\matstat(\sigma)=3$ and $\hat{\imath}(\sigma)=1$.

\begin{defi}
    A \emph{$31-2$ pattern} in $\sigma \in \mathfrak S_n$ is a triple $1\leq i < i+1 < j \leq n$ such that $\sigma_j> \sigma_i > \sigma_{i+1}$. We denote by  $31-2(\sigma)$ the number of $31-2$ patterns in $\sigma$.
\end{defi}
In \cite{HanRandrianarivonyZeng1999}, the authors introduced an interesting $t$-analogue to the Euler numbers that was subsequently studied in \cite{Chebikin2008} and \cite{JosuatVerges2010}. 

%TODO power of t missing
\begin{defi}\label{def:qeuler}
    For all $n\in \N$, define:
    \begin{equation}
        E_n(t) \coloneqq \sum_{\sigma\in \mathfrak A_n} t^{31-2(\sigma)}.
    \end{equation}
\end{defi}

This polynomial $E_n(t)$ has several beautiful properties including the facts that the generating functions
$\sum_{n\ge 0} E_{2n}(t)z^n$
and  $\sum_{n\ge 0} E_{2n+1}(t)z^n$ have nice continued fraction expressions \cites{HanRandrianarivonyZeng1999, 
JosuatVerges2010}.

Here we study a shift of this $t$-analogue, namely $t^{\lfloor n^2/4\rfloor}E_n(t)$.
This $t$-analogue is naturally connected to our 3-inversion statistics.
\begin{prop}\label{prop:t-euler-alternating-matstat}
For all $n\in \N$ 
    \begin{equation}
        t^{\lfloor n^2/4\rfloor}E_n(t) = \sum_{\sigma\in \mathfrak A_n} t^{\matstat(\sigma)}.
    \end{equation}
\end{prop} 

\begin{proof} 
This is a Corollary of Lemma 9.4 of \cite{Chebikin2008}. In this lemma, Chebikin proves that if $\sigma\in \mathfrak A_n$ then $\lfloor n^2/4\rfloor+31-2(\sigma)=\hat{\imath}(\sigma)$ and we just remarked that if  $\sigma\in \mathfrak A_n$
then $\matstat(\sigma)=\hat{\imath}(\sigma)$.
\end{proof}

\section{Generating trees}\label{sec:trees}

The goal of this section is to prove the following. Recall that we denote by $\mathfrak S^{\bullet}_n(1^n)$ the subset of $\mathfrak{S}_n^{\bullet}$ of marked permutations with schedule $1^n$. 

\begin{thm}\label{thm:trees-isom}
	There is a bijection $\phi : \mathfrak{S}_n^{\bullet}(1^n) \to \mathfrak{S}_n$ with the following properties: 
  \begin{enumerate}[(i)]
    \item $\revmaj(\tau) = \matstat(\phi(\tau))$,
    \item $\mathsf{dec}(\tau) = \monot(\phi(\tau))$. 
  \end{enumerate}
\end{thm}

We will exhibit this bijection by constructing two isomorphic generating trees. 

\subsection{The tree for decorated permutations of schedule \texorpdfstring{$1^n$}{1n}}
We continue to use Convention~\ref{conv:preceding-zero}: for any $\sigma\in \mathfrak{S}_n$, we set $\sigma_0=0$. 

Let $\tau\in{\mathfrak S}_n$. Using Definition~\ref{def:cyclic-runs}, we know that $\tau_0\ldots \tau_j$ is a cyclic run if and only if $\tau_1>\ldots >\tau_j$. The generation of the tree will rely on the following manipulation. 
\begin{defi}
    Let $\tau\in{\mathfrak S}_n$ and $1\leq k \leq n$. Define $k + \tau$ to be the permutation of $\mathfrak{S}_n$ defined by 
    \begin{align*}
        \begin{cases}
            \sigma_0 = 0 \\
            \sigma_{i+1}=(k+\tau_{i})\mod n+1\quad {\rm for}\ 0\le i\le n.
        \end{cases}
    \end{align*} 
\end{defi}
\begin{example*}
    % decorate 4 and 5 for schedule 1^n 
    If $n=7$, $\tau = 01423657$, and $k = 3$, we have 
    \[
        k + \tau = 034756182.
    \]
\end{example*}

\begin{lemma}
    Let $\tau\in\mathfrak S_n$ and $\sigma\coloneqq k + \tau$. For $0\le i\le j\le n$, $\tau_i\cdots \tau_j$ is a cyclic run of $\tau$ if and only if  $\sigma_{i+1}\cdots \sigma_{j+1}$ is a cyclic run of $\sigma$. Moreover, for $1\le j\le n$, ($\tau_0\cdots \tau_j$ is a cyclic run and $\tau_j+ k> n+1$) if and only if  $\sigma_{0}\cdots \sigma_{j+1}$ is a cyclic run.
    \label{lem:cyclic}
\end{lemma}
\begin{proof}
    The first statement is a direct consequence of the definition of a cyclic run (Definition~\ref{def:cyclic-runs}).
    Note that for any permutation $\tau$, $\tau_{0}\ldots \tau_{j}$ is a cyclic run if and only if  $\tau_1>\ldots >\tau_{j}$. Therefore $\sigma_{0}\ldots \sigma_{j+1}$ is a cyclic run
    if and only if $k=\sigma_1>\ldots >\sigma_{j+1}$. As $\sigma_{j+1}=(\tau_j +k)\mod n+1$, we have that $\sigma_{j+1}<k$ if and only if $\tau_j +k>n+1$.
\end{proof}

The structure of the generating tree will closely depend upon the following quantity.
\begin{defi}\label{def:struct-att-tree1}
    Given $\tau\in \mathfrak{S}_n^{\bullet}(1^n)$, we define its \emph{structural attribute}, $a(\tau)$ to be the value of its first undecorated letter.
\end{defi}

\begin{defi}\label{def:tree1}
    Define a tree $\mathcal T_1$ of decorated permutations as follows. Take its root to be $1\in \mathfrak S_1$. 
    For $\tau\in \mathfrak{S}_n^{\bullet}(1^n)$ with $n\geq 1$  and $1\leq k \leq n+1$, the \emph{$k$-th descendant of  $\tau$ in $\mathcal{T}_1$}, denoted $\delta_k(\tau)$ is the decorated permutation whose underlying permutation is $\sigma = k+\tau$. The decorations of $\sigma$ are as follows: 
    \begin{itemize}
        \item $\sigma_0 = 0$ is never decorated
        \item $\sigma_1$ is decorated if and only if $k>n+1-a(\tau)$
        \item for $i\geq 2$, $\sigma_i$ is decorated if and only if $\tau_{i-1}$ is decorated.
    \end{itemize} 
\end{defi}
See the tree on the left in Figure~\ref{fig:trees} for an illustration of the first three levels of this tree.
\begin{example}
    For example, take $n=5$ and $\tau = 0\dt 534\dt 21$. Then $a(\tau) = 3$ and so $n+1-a(\tau) = 3$. We have 
    \[\delta_2(\tau) ={0}2\dt 156\dt 43 \hspace{2cm} \delta_4(\tau) = {0}\dt 4 \dt 312 \dt 65.\]
\end{example}

% \begin{defi}
%     For $\tau\in \mathfrak{S}_n^{\bullet}(1^n)$ and $1\leq k \leq n+1$ we define $(k+\tau  \mod n+1)$ to be the marked word whose $i$-th letter is the unique integer in $\{1,\dots,n+1\}$ congruent to $ k + \tau_i \mod n+1$; and to be decorated if and only if $\tau_i$ is decorated.  
% \end{defi}

Let us prove some observations about $\mathcal T_1$. 
Suppose that $\tau$ is a decorated permutation and its schedule is $1^{n}$. 
\begin{lemma}
Take $\tau\in \mathfrak{S}_n^{\bullet}(1^n)$. If $\tau_{0}\ldots \tau_{j}$ is a right-maximal cyclic run then there is a unique $\ell$ such that $1\le \ell\le j$ and $\tau_{\ell}$ is undecorated. 
\label{lem:unique-undecorated}
\end{lemma}
\begin{proof}
    Call $r \coloneqq \tau_{0}\ldots \tau_{j}$. If $r$ it contained at least two undecorated letters, the second largest one would have schedule $2$. If $r$ contains no undecorated letters there must be exactly one undecorated $\tau_i$ with $i>j$ in the maximal decreasing cyclic run starting at $\tau_k$ for all $1\leq k\leq j$.  It follows that the maximal decreasing run ending at $\tau_i$ contains only decorated letters and so its schedule would be $0$ (see Remark~\ref{rem:cyclic-runs}). 
\end{proof}

\begin{prop}
    The $n$-th level of the tree $\mathcal{T}_1$ contains exactly all the decorated permutations whose schedule numbers are $1^n$. 
\end{prop}
\begin{proof}
    Take $\tau\in \mathfrak S_n^\bullet(1^n)$. We have to show that for all $k\in [n+1]$, $\sched(\delta_k(\tau)) = 1^{n+1}$. The result will then follow from Lemma~\ref{lem:unique-decoration-1n} and the fact that the tree in Definition~\ref{def:tree1} clearly generates all permutations of $n$ at the $n$-th level. 

    Set $\tau_\ell\coloneqq a(\tau)$, the unique undecorated entry in the  right-maximal cyclic run containing $\tau_0$. Thanks to Lemma \ref{lem:unique-undecorated}, this $\ell$ is unique.

    We construct the decorated permutation $\sigma=\delta_k(\tau)$.

    Thanks to Lemma \ref{lem:cyclic} we know that if $\tau_0\ldots \tau_j$ is not a cyclic run then the cyclic runs that contain $\sigma_{j+1}$ are in bijection with the cyclic runs that contain $\tau_j$ and the schedule of $\sigma_{j+1}$ is equal to the schedule of $\tau_j$, which is one.

    Moreover if $j>0$ and $\tau_0\ldots \tau_j$ is a cyclic run and $\tau_j$ is decorated, the cyclic runs that start at $\sigma_{j+1}$ are in bijection with the cyclic runs that start at $\tau_j$ and the schedule of $\sigma_{j+1}$ is equal to the schedule of $\tau_j$, which is one.

    We just have to show that $\sigma_{\ell+1}$ and $\sigma_1$ have schedule 1.

    If $\tau_\ell+k>n+1$ then (by Lemma~\ref{lem:cyclic}), $\sigma_{0}\ldots \sigma_{\ell+1}$ is a cyclic run.
    As
    $\sigma_{\ell+1}$ is not decorated, we must decorate $\sigma_1$ to force that $\sigma_{\ell+1}$ and $\sigma_1$ have schedule 1.

    If $\tau_\ell+k\le n+1$ then $\sigma_{0}\ldots \sigma_{\ell+1}$ is not a cyclic run but $\sigma_{1}\ldots \sigma_{\ell+1}$ is. Therefore  $\sigma_{1}\ldots \sigma_{\ell+1}$ is a left maximal run and to force $\sigma_{\ell+1}$ to have schedule 1, $\sigma_1$ must be non-decorated (otherwise $\sigma_{\ell+1}$ has schedule 0). Moreover if $\sigma_1$ is not decorated, its schedule is one and we can conclude the proof.
    % In general, if $\tau \in \mathfrak S_n^\bullet$ and the first decreasing run is of length $r$, we have \[\sched(\tau)_i = \sched((k+\tau) \mod n+1)_i \quad \text{if $i>r$ or $\tau_i$ is decorated.}\]  This is clear from the definition of a schedule number in terms of cyclic runs (Definition~\ref{rem:cyclic-runs}), which is invariant by the operation $\tau\mapsto (k+\tau) \mod n+1$, provided that the schedule we want to compute is not an undecorated letter in the first run (this corresponds to a rotation of Figure~\ref{fig:cyclic-runs}, viewed as a cylinder, as explained in Example~\ref{ex:schedules}).  
    % Now we take $\tau\in \mathfrak S_n^\bullet$ with $\sched(s) = 1^n$. Since all the schedules are $1$, it follows from the schedule Definition~\ref{def:schedule} that the first decreasing run of $\tau$ contains exactly one undecorated letter. Indeed if it contained at least two, the second largest one would have schedule $2$; and if it contained none, there must be a second run and its first letter would have schedule $0$. Call $a(\tau)$ this only undecorated letter in the first run of $\tau$. Similarly, for $\delta_k(\tau)$ to have schedule $1^{n+1}$ the first run must contain exactly one undecorated letter. So the first letter of $\delta_k(\tau)$ may only be undecorated if $k <(a(\tau) + k) \mod n+1$, in other words, if $k\leq n+1-a(\tau)$. For $k>n+1-a(\tau)$, the first letter of $\delta_k(\tau)$ is decorated and it is clear that both $\dt k$ and $a(\tau)$ have schedule $1$ in this case.  
\end{proof}

\begin{prop}\label{prop:revmaj}
   For all decorated permutations $\tau$, we have \[\revmaj(\delta_k(\tau)) = \revmaj(\tau) + (n+1-k).\]
\end{prop}

\begin{proof}
	This is clear for $k=n+1$, as $\delta_{n+1}(\tau)$ is obtained by appending $n+1$ at the beginning of $\tau$ (which creates no new descent in the reversed permutation, thus preserving the reverse major index). 
	
    Next, we show $\revmaj(\delta_{k}(\tau)) =\revmaj(\delta_{k-1}(\tau)) -1$. To go from $\delta_{k-1}(\tau)$ to $\delta_{k}(\tau)$, we add $1$ to every letter, and replace $n+1$ with $1$. Let $m$ be the index of $n+1$ in $\delta_{k-1}(\tau)^\rev$. Then we have:
    \[\Des(\delta_{k}(\tau)^\rev) = 
    \begin{cases}
        \left(\Des(\delta_{k-1}(\tau)^\rev) \setminus \{m\}\right) \cup \{m-1\} &\text{if $m> 1$}, \\
        \Des(\delta_{k-1}(\tau)^\rev) \setminus \{m\} &\text{if $m=1$}.
    \end{cases}
    \] 
    Indeed, the $m$-th letter of $\delta_{k-1}(\tau)^\rev$ was equal to $n+1$ and so was a descent, but in $\delta_{k}(\tau)^\rev$ it equals $1$ and so is not descent. Furthermore if $m>1$, the $m-1$-th letter of $\delta_{k}(\tau)^\rev$ is followed by $1$ and so must be a descent. In any case 
    \begin{align*}
        \revmaj(\delta_{k}(\tau)^\rev) = \revmaj(\delta_{k-1}(\tau)^\rev) -m + m -1 = \revmaj(\delta_{k-1}(\tau)^\rev) -1
    \end{align*} and the result follows. 
\end{proof}

Finally, the following property follows easily from Definitions~\ref{def:struct-att-tree1} and~\ref{def:tree1}.

\begin{prop}\label{prop:struct-tree1}
    Recall that $\mathsf{dec}(\tau)$ is the number of decorations of $\tau\in \mathfrak{S}_{n}^{\bullet}(1^n)$.  We have: 
    \[
        \mathsf{dec}(\delta_k(\tau)) = \mathsf{dec}(\tau) + \chi(n+1-k\leq a(\tau)).
    \]
    Furthermore, we have:
    \[
        a(\delta_k(\tau))= 
        \begin{cases}
            k &\text{if } k\leq n+1 - a(\tau), \\
            a(\tau) &\text{if } k> n+1- a(\tau).
        \end{cases}
    \]
\end{prop}

\subsection{A second tree related to peaks and valleys}

\begin{defi}\label{def:tree2}
    Define a tree $\mathcal T_2$ of permutations as follows. Take its root to be $1\in \mathfrak{S}_1$. 
    For $\sigma\in \mathfrak{S}_n$ and $1\leq l \leq n+1$, define $\eta_l(\sigma) \in \mathfrak{S}_{n+1}$ to be the unique permutation $\sigma'$ such that $\sigma'_{n+1}=l$ and $\sigma'_1, \dots,\sigma'_n$ are in the same relative order as $\sigma_1, \dots,\sigma_n$.  In other words 
    \[\sigma'_i = 
   \sigma_i + \chi(\sigma_{i} \geq l); 
    \] for $1\leq l \leq n+1$.
    This permutation is \emph{called the insertion of $l$ in $\sigma$}. The $\eta_l(\sigma)$ will form the descendants of $\sigma$ in $\mathcal T_2$.   
\end{defi}
See the tree on the right in Figure~\ref{fig:trees}. Later (Definition~\ref{def:kth-des-tree2}), we will define a total order on the descendants of a node in $\mathcal T_2$, in a way that will give the isomorphism with $\mathcal T_1$. This ordering will closely depend upon the following quantity. 

\begin{defi}
    Given $\sigma \in \mathfrak S_n$, define it \emph{structural attribute}:
    \[\tilde{a}(\sigma) = 
    \begin{cases}
        n+1 -\sigma_n & \text{if } \sigma_{n-1} <\sigma_n,\\
        \sigma_n & \text{if } \sigma_{n-1} >\sigma_n,
    \end{cases},\] 
    where we consider $\sigma_0=0$ in case $n=1$.
\end{defi}

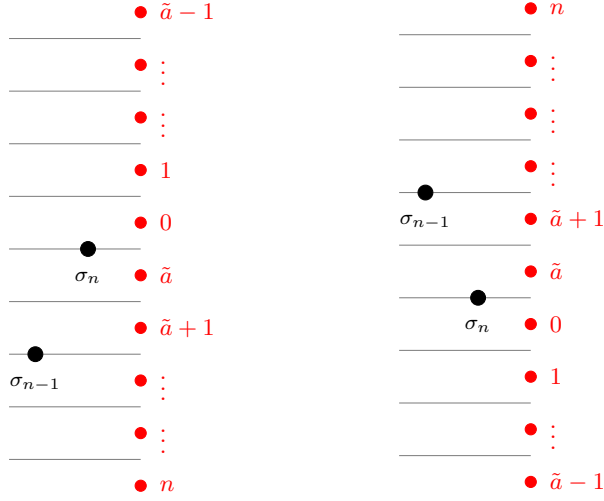
\begin{figure}
    \centering
    \begin{tikzpicture}[scale=.7]
        \foreach \i in {1,...,9}{ 
        \draw[opacity=.5] (-.5,\i) -- (2,\i);
        }     
        \foreach \j in {1, ..., 10}{
            \draw[red,fill] (2,\j - .5) circle (3pt);
            \node (\j) at (2.2, \j -.5) {};   
        }
        \textcolor{red}{\small
        \node[right] at (6) {$0$};
        \node[right] at (7) {$1$};
        \node[right] at (8) {$\vdots$};
        \node[right] at (9) {$\vdots$};
        \node[right] at (10) {$\tilde{a}-1$};
        \node[right] at (5) {$\tilde{a}$};
        \node[right] at (4) {$\tilde{a}+1$};
        \node[right] at (3) {$\vdots$};
        \node[right] at (2) {$\vdots$};
        \node[right] at (1) {$n$};
        }
        \draw[fill] (0,3) circle (4pt) node[below = 5pt]{\footnotesize$\sigma_{n-1}$};
        \draw[fill] (1,5) circle (4pt) node[below = 5pt]{\footnotesize$\sigma_{n}$};
        %\draw[decorate, decoration={brace, amplitude = 8pt, mirror}] ($(6) + (1.3,-.3)$) -- ($(10) + (1.3,.3)$) node[pos=.5,right=.3cm]{$\tilde{a} = n+1-\sigma_{n}$};
    \end{tikzpicture}
    \hspace{2cm}
    \begin{tikzpicture}[scale=.7]
        \foreach \i in {1,...,9}{ 
        \draw[opacity=.5] (-.5,\i) -- (2,\i);
        }     
        \foreach \j in {1, ..., 10}{
            \draw[red,fill] (2,\j - .5) circle (3pt);
            \node (\j) at (2.2, \j -.5) {};   
        }
        \textcolor{red}{\small
        \node[right] at (4) {$0$};
        \node[right] at (3) {$1$};
        \node[right] at (2) {$\vdots$};
        \node[right] at (1) {$\tilde{a}-1$};
        \node[right] at (5) {$\tilde{a}$};
        \node[right] at (6) {$\tilde{a}+1$};
        \node[right] at (7) {$\vdots$};
        \node[right] at (8) {$\vdots$};
        \node[right] at (9) {$\vdots$};
        \node[right] at (10) {$n$};
        }
        \draw[fill] (0,6) circle (4pt) node[below = 5pt]{\footnotesize$\sigma_{n-1}$};
        \draw[fill] (1,4) circle (4pt) node[below = 5pt]{\footnotesize$\sigma_{n}$};
        %\draw[decorate, decoration={brace, amplitude = 8pt, mirror}] ($(1) + (1.3,-.3)$) -- ($(4) + (1.3,.3)$) node[pos=.5,right=.3cm]{$\tilde{a}=\sigma_{n}$};
    \end{tikzpicture}
    \caption{Contribution to $\matstat$ for the $n+1$ different insertions into a permutation of $n$.}\label{fig:insertionT2}
\end{figure}

\begin{prop}\label{prop:matstat-insertion}
%    For $\sigma\in \mathfrak S_n$ there exists a bijection \[\psi : \{1,2,\dots, n+1\}\rightarrow \{0,1,\dots, n\}\] such that $\matstat(\eta_l(\sigma)) = \matstat(\sigma) + \psi(l)$. 
    For $\sigma\in \mathfrak S_n$, the map $\psi$ on $\{1,2,\dots, n+1\}$ defined by \[\matstat(\eta_l(\sigma)) = \matstat(\sigma) + \psi(l)\] is a bijection onto $\{0,1,\dots, n\}$.
\end{prop}
\begin{proof}
    Take $\sigma \in \mathfrak S_n$ and set $\sigma' = \eta_l(\sigma)$. Notice that for $j\leq n$ we have that $(i,j)$ forms a $3$-inversion in $\sigma$ if and only if it forms a $3$ inversion in $\sigma'$ (see Definitions~\ref{def:3inversion} and \ref{def:tree2}). It follows that $\matstat(\sigma')=\matstat(\sigma') + s$ where $s$ is the number of $3$-inversions of the form $(i,n+1)$.   
    
    We distinguish two cases. We have drawn a schematic representation of the proof in Figure~\ref{fig:insertionT2}. 
    \begin{enumerate}
        \item \textit{Either $\sigma$ ends with an ascent}, that is $\sigma_{n-1}<\sigma_n$. So we have $\tilde a \coloneqq \tilde{a}(\sigma) = n+1-\sigma_n$.  This is illustrated in the left part of Figure~\ref{fig:insertionT2}, where the values in red are the value of $\psi(l)$.
        \begin{itemize}
            \item For $l>\sigma_n $, $\sigma'_{n+1}$ is a double ascent. The number of $i$'s such that $\sigma'_{n} <\sigma'_i <\sigma'_{n+1}$ is equal to $l-(\sigma_n +1)$. So for $l=\sigma_n +1,\sigma_n + 2, \dots, n+1$ we have $\psi(l) = 0, 1, \dots, \tilde{a}-1$, respectively.
            \item For $l\leq\sigma_n$, $\sigma'_{n+1}$ is a peak. The number of $i$'s such that $\sigma'_{i} >\sigma'_{n+1}$ is equal to $n + l - 1$. So for $l =\sigma_n,\sigma_n-1,\dots, 1$ we have $\psi(l) = \tilde{a},\tilde{a}+1, \dots, n$, respectively. 
        \end{itemize}
        \item \textit{Or, $\sigma$ ends with a descent}, that is $\sigma_{n-1}>\sigma_n$. So we have $\tilde a \coloneqq \tilde{a}(\sigma)=\sigma_n$.  This is illustrated in the right part of Figure~\ref{fig:insertionT2}.
        \begin{itemize}
            \item For $l\leq\sigma_{n}$, $\sigma'_{n+1}$ is a double descent. The number of $i$ such that $\sigma'_{n} >\sigma'_{i} >\sigma'_{n+1}$ equals $\sigma_n - l$. So for $l =\sigma_n,\sigma_{n}-1, \dots, 1$ we have $\psi(l) = 0,1,\dots, \tilde{a}-1$, respectively. 
            \item For $l>\sigma_n$, $\sigma'_{n+1}$ is a valley. The number of $i$ such that $\sigma'_{i} >\sigma'_{n+1}$ is equal to $l-1$. So for $l=\sigma_n + 1, \dots, n+1$, we have $\psi(l) = \tilde{a}, \tilde{a}+1, \dots, n$, respectively. 
        \end{itemize}
    \end{enumerate}
In each of the two cases, we see that the values taken by $\psi(l)$ are exactly $0,\dots,n$.
\end{proof}

\begin{defi}\label{def:kth-des-tree2}
    For $\sigma\in\mathfrak S_n$ we define the \emph{$k$-th descendant of  $\sigma$ in $\mathcal{T}_2$}, $\tilde \delta_k(\sigma)$, to be the unique descendant $\eta_l(\sigma)$ such that $\psi(l) = n+1-k$, in other words $\matstat(\tilde\delta_k(\sigma)) = \matstat(\sigma) + (n+1-k)$.
\end{defi}

Finally, we can deduce the following from the proof of Proposition~\ref{prop:matstat-insertion}.
\begin{prop}\label{prop:struct-tree2}
    Recall that $\mathsf{monot}(\sigma)$ denotes the number of double ascents and descents of $\sigma\in \mathfrak{S}_{n}$.  We have:
    \[
        \mathsf{monot}(\eta_k(\sigma)) = \mathsf{monot}(\sigma) + \chi(n+1-k\leq \tilde a(\sigma)).
    \]
    Furthermore, we have: 
    \[
        \tilde a(\eta_k(\sigma))= 
        \begin{cases}
            k &\text{if } k\leq n+1 - \tilde a(\sigma), \\
            \tilde a(\sigma) &\text{if } k> n+1- \tilde a(\sigma).
        \end{cases}
    \]
\end{prop}

\subsection{Isomorphism of the trees}

To prove Theorem~\ref{thm:trees-isom}, we define $\phi$ recursively as follows:
\begin{align*}
    &\phi(1) = 1,\\
    &\phi(\delta_k(\sigma)) = \eta_k(\phi(\sigma)) \text{ for all } \sigma\in\biguplus_{n\geq 1} \mathfrak S_n^{\bullet}(1^n).
\end{align*}
In Figure~\ref{fig:trees}, the $k$-th descendent of each
node is the $k$-th descendant from the bottom, so that the image by $\phi$ of an element in the left tree can be obtained by looking for the corresponding element in the right tree, if we were to ``superpose'' one tree on the other. 

From this definition of $\phi$ and the definition of \emph{$k$-th descendant} as being the descendant adding $n+1-k$ units to the relevant statistic in both trees (See Definition~\ref{def:tree1}, Proposition~\ref{prop:revmaj} and Definition~\ref{def:kth-des-tree2}), we may conclude that \[\revmaj(\sigma) = \matstat(\phi(\sigma)).\] 
From Proposition~\ref{prop:struct-tree1} and Proposition~\ref{prop:struct-tree2}, we deduce \[\mathsf{dec}(\sigma) = \mathsf{monot}(\phi(\sigma)).\]
Thus, we have now established Theorem~\ref{thm:trees-isom}.

\tikzset{level 1/.style = {sibling distance = 9.4cm},
level 2/.style = {sibling distance = 3.15cm},
level 3/.style = {sibling distance = .75cm, level distance = 2cm},
grow = right}
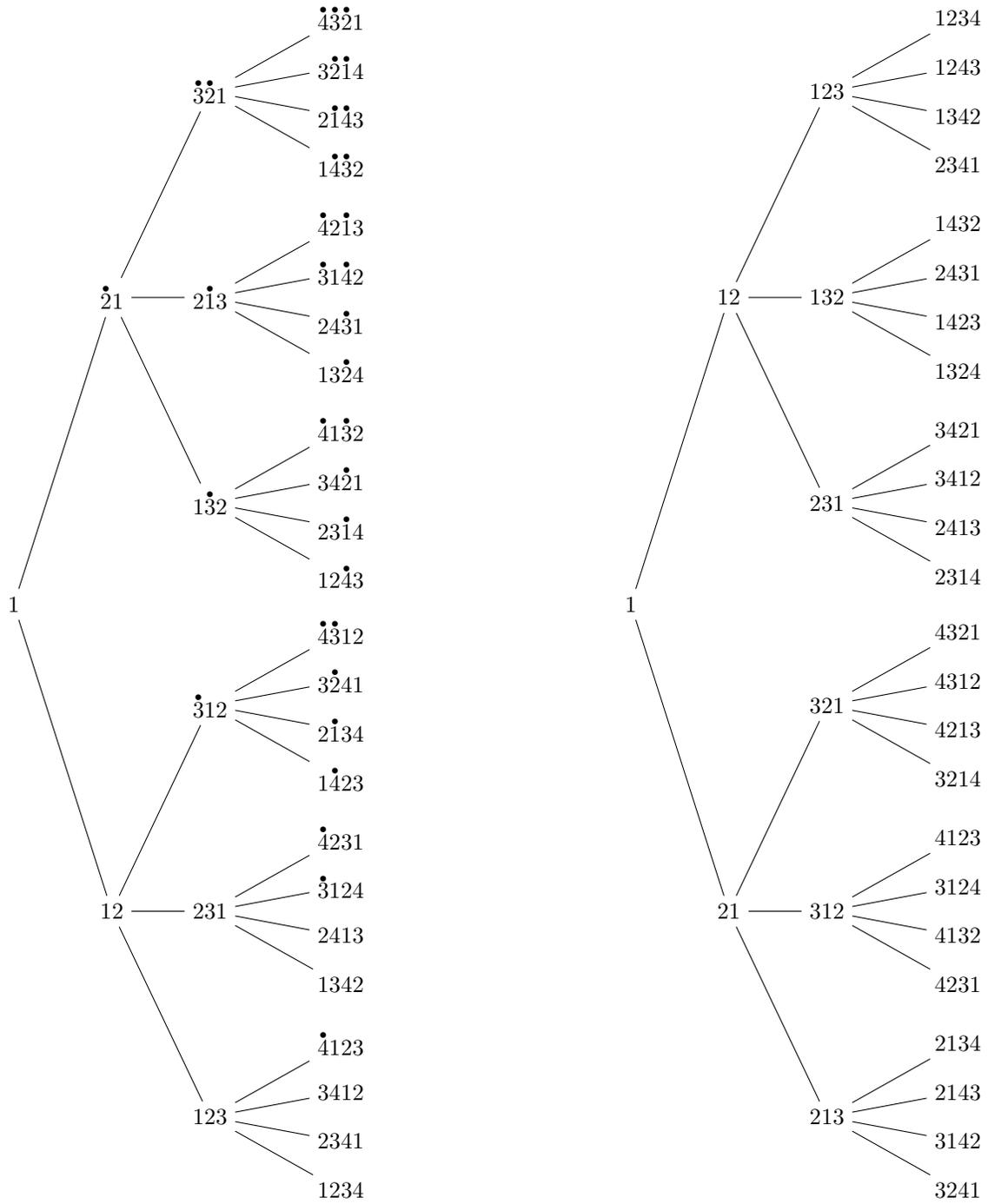
\begin{figure}
    \centering
    \begin{tikzpicture}
        \node {$1$}
        child {node {$12$} 
            child {node{$123$}
                child {node {$1234$}}
                child {node {$2341$}}
                child {node {$3412$}}
                child {node {$\dt 4123$}}
                }
            child {node{$231$}
                child {node {$1342$}}
                child {node {$2413$}}
                child {node {$\dt 3124$}}
                child {node {$\dt 4231$}}
            }
            child {node{$\dt 312$}
                child {node {$1\dt 423$}}
                child {node {$2\dt 134$}}
                child {node {$3\dt 241$}}
                child {node {$\dt 4\dt 3 12$}}
            }
        }
        child {node {$\dt 2 1$}
            child {node{$1\dt 32$}
                child {node {$12\dt 43$}}
                child {node {$23\dt 14$}}
                child {node {$34\dt 21$}}
                child {node {$\dt 41 \dt 32$}}
                }
            child {node{$2\dt 13$}
                child {node {$13\dt 2 4$}}
                child {node {$24\dt 3 1$}}
                child {node {$\dt 3 1 \dt 4 2$}}
                child {node {$\dt 4 2 \dt 1 3$}}
            }
            child {node{$\dt 3\dt 2 1$}
                child {node {$1\dt 4 \dt 3 2$}}
                child {node {$2 \dt 1 \dt 4 3$}}
                child {node {$3 \dt 2 \dt 1 4$}}
                child {node {$\dt 4 \dt 3 \dt 2 1$}}
            }
        };
    \end{tikzpicture}
    \hfill
    \begin{tikzpicture}
        \node {$1$}
        child {node {$21$}
            child {node{$213$}
                child {node {$3241$}}
                child {node {$3142$}}
                child {node {$2143$}}
                child {node {$2134$}}
            }
            child {node{$312$}
                child {node {$4231$}}
                child {node {$4132$}}
                child {node {$3124$}}
                child {node {$4123$}}
            }
            child {node{$321$}
                child {node {$3214$}}
                child {node {$4213$}}
                child {node {$4312$}}
                child {node {$4321$}}
            }
        }
        child {node {$12$}
            child {node{$231$}
                child {node {$2314$}}
                child {node {$2413$}}
                child {node {$3412$}}
                child {node {$3421$}}
                }
            child {node{$132$}
                child {node {$1324$}}
                child {node {$1423$}}
                child {node {$2431$}}
                child {node {$1432$}}
            }
            child {node{$123$}
                child {node {$2341$}}
                child {node {$1342$}}
                child {node {$1243$}}
                child {node {$1234$}}
            }
        };
    \end{tikzpicture}
    \caption{$\mathcal{T}_1$ and $\mathcal{T}_2$ up to level 4. }\label{fig:trees}
\end{figure}

\section{Conclusion and future directions}

Let us now put the pieces together.
From the bijection in Theorem~\ref{thm:trees-isom}, we may deduce that
\[
    \sum_{\tau\in\mathfrak S_n^{\bullet }(1^n)} t^{\revmaj(\tau)}z^{\mathsf{dec}(\sigma)} = \sum_{\sigma\in \mathfrak{S}_n} t^{\matstat(\sigma)}z^{\monot(\sigma)}.
\]
Recall Equation~\ref{eq:combinatorics-Delta-q-1}:
\[\sum_{P\in \stLD(n)^{\bullet k}} (-1)^{\dinv(P)}t^{\area(P)}= \sum_{\tau\in\mathfrak S_n^{\bullet k}(1^n)} t^{\revmaj(\tau)}.\]
Combining these last two equations and summing over the $k$ we get: 
\[\sum_{k = 0}^{n-1}\left(\sum_{P\in \stLD(n)^{\bullet k}}t^{\area(P)}(-1)^{\dinv(P)} \right)z^k= \sum_{\sigma \in \mathfrak S_n} t^{\matstat(\sigma)}z^{\monot(\sigma)};\]
which is exactly the statement in Theorem~\ref{thm:main}.

Since the left hand side of \eqref{eq:combinatorics-Delta-q-1} is exactly the combinatorics of the Hilbert series of the Delta conjecture (Equation~\ref{eq:hilbert-delta}) at $q=-1$, the truth of the Delta conjecture would imply 
\begin{align*}
    \sum_{k=0}^{n-1}\left. \langle \Delta'_{e_{n-k-1}}e_n,h_1^n\rangle\right|_{q=-1}z^k = \sum_{\sigma\in\mathfrak S_n} t^{\matstat(\sigma)}z^{\monot(\sigma)},
\end{align*}
as announced in Corollary~\ref{cor:main}. 

Since at $k=0$ the Delta conjecture reduces to the shuffle theorem, our result implies that 
\[
    \langle \nabla e_n, h_1^n\rangle |_{q=-1}
        = 
    \sum_{\substack{\sigma\in\mathfrak S_n \\ \monot(\sigma) = 0}} t^{\matstat(\sigma)} = \sum_{\sigma\in \mathfrak{A}_n} t^{\matstat(\sigma)}= t^{\lfloor n^2/4\rfloor}E_n(t),
\] 
where the last equality comes from Proposition~\ref{prop:t-euler-alternating-matstat}. Thus we have established \eqref{eq:t-euler}.

Computational evidence suggests that the evaluation at $q=-1$ yields $t$-positive results for many other polynomials related to the shuffle theorem and Delta conjecture. For example 
\begin{itemize}
    \item $\langle \nabla e_n, h_\mu\rangle$ for any partition $\mu$;
    \item $\langle \nabla \omega(p_n);h_1^n\rangle$, where $\nabla \omega( p_n)$ is the symmetric function related to the square theorem (conjectured in \cite{LoehrWarrington2007} and proved in \cite{Sergel2017});
    \item $\langle \nabla E_{n,k}, h_1^n\rangle$ where $E_{n,k}$ is the symmetric function refinement of $e_n$ introduced in \cite{GarsiaHaglund2002};
    %\item $\langle \Theta_{e_k}\Theta_{e_l} \nabla e_{n-k-l}, h_1^n\rangle$, where $\Theta_{e_k}\Theta_{e_l} \nabla e_{n-k-l}$ is the symmetric function  introduced in \cite{DAdderioIraciVandenWyngaerd2021} that conjecturally unifies both versions of the Delta conjecture;
    %\item $\langle \Theta_{e_\lambda} e_1, h_1^{\vert \lambda \vert +1}\rangle$, see the Theta trees conjecture in \cite{D’AdderioIraciLeBorgneRomeroVandenWyngaerd2022}. 
\end{itemize}

It would thus be interesting to study the $q=-1$ evaluation in a more general framework; for example in modified Macdonald polynomials.

Let us end this section with a combinatorial problem.  Consider the sums
\begin{equation} \label{def_dnj}
    D_{n,j}(t) 
        = 
    \sum_{\substack{\sigma \in \mathfrak{S}_n, \\ \monot(\sigma) \leq j} } t^{\matstat(\sigma)}.
\end{equation}
In particular, this is $0$ if $j<0$, and $D_{n,j}(q) = D_{n,n-1}(q)$ if $j>n-1$.  

\begin{conj}
    For $n>1$, we have:
    \begin{equation}  \label{rec_deux_termes}
        D_{n,j}(t) = t^{n-j-1} [j+1]_t D_{n-1,j+1}(t) + [n-j-1]_t D_{n-1,j-1}(t).
    \end{equation}
\end{conj}

It would be very interesting to give a combinatorial proof of~\eqref{rec_deux_termes}.  Of course, one might try do to this starting from the other combinatorial interpretation: 
\[
    D_{n,j}(q) = \sum_{\substack{ \sigma \in \mathfrak{S}^\bullet(1^n), \\ \mathsf{dec}(\sigma)\leq j }} q^{\revmaj(\sigma)}.
\]
In the case $t=1$, it is possible to give a combinatorial proof of the conjecture, based on the combinatorial interpretation in~\eqref{def_dnj}.  The idea is to consider the map $\sigma \mapsto \sigma'$ (where, for $\sigma\in\mathfrak{S}_n$, $\sigma'\in\mathfrak{S}_{n-1}$ is obtained by removing the entry $n$), and examine how $\monot$ is distributed among the $n$ pre-images of a given $\sigma \in\mathfrak{S}_{n-1}$.

\section*{Acknowledgements}
The authors would like to acknowledge the support of ANR Combiné ANR-19-CE48-0011. 
The last author would like to thank the FSMP, who supported her during this research.
Thank you to Michele D'Adderio for insightful discussions. 

\bibliographystyle{amsalpha}
\bibliography{q-euler.bib}

\end{document}